\documentclass[11pt]{amsart}
\usepackage{
amssymb
}
\newcommand{\no}[1]{#1}

\renewcommand{\no}[1]{}  \newcommand{\upDelta}{\Delta} 

\no{\usepackage{times}\usepackage[
subscriptcorrection, slantedGreek, 
nofontinfo]{mtpro}
\renewcommand{\Delta}{\upDelta}
}

\newtheorem{theorem}{Theorem}
\newtheorem{proposition}{Proposition}
\newtheorem{lemma}{Lemma}
\newtheorem{definition}{Definition}
\newtheorem{corollary}{Corollary}

\theoremstyle{remark}
\newtheorem{remark}{Remark}

\DeclareMathOperator{\supp}{supp}

\newcommand{\eps}{\varepsilon}
\newcommand{\R}{{\bf R}}
\newcommand{\Id}{\mbox{Id}}
\renewcommand{\r}[1]{(\ref{#1})}
\newcommand{\PDO}{$\Psi$DO}
\newcommand{\be}[1]{\begin{equation}\label{#1}}
\newcommand{\ee}{\end{equation}}

\DeclareMathOperator{\n}{neigh}
\renewcommand{\d}{\mathrm{d}}

\newcommand{\bo}{\partial M}
\newcommand{\Mint}{M^\text{\rm int}}
\newcommand{\D}{\mathcal{D}}

\title[Lens rigidity]{Local lens rigidity with incomplete data for a class of\\ non-simple Riemannian manifolds}

\author[P. Stefanov]{Plamen Stefanov}
\address{Department of Mathematics, Purdue University, West Lafayette, IN 47907}
\thanks{First author partly supported by NSF Grant DMS-0400869}
\author[G. Uhlmann]{Gunther Uhlmann}
\address{Department of Mathematics, University of Washington, Seattle, WA 98195}
\thanks{Second author partly supported by NSF and a Walker Family Endowed Professorship}

\begin{document}

\begin{abstract} 
Let $\sigma$ be the scattering relation on a compact   Riemannian manifold $M$ with non-necessarily convex boundary, that maps initial points of geodesic rays on the boundary and initial directions to the outgoing point on the boundary and the outgoing direction. Let $\ell$ be the length of that geodesic ray. We study the question of whether the metric $g$ is uniquely determined, up to an isometry, by knowledge of $\sigma$ and $\ell$ restricted on some subset $\D$. We allow possible conjugate points but we assume that the conormal bundle of the geodesics issued from $\D$ covers $T^*M$; and that those geodesics have no conjugate points. Under an additional topological assumption, we prove that $\sigma$ and $\ell$ restricted to $\D$ uniquely recover an isometric copy of  $g$ locally near  generic metrics, and in particular, near real analytic ones.
\end{abstract}

\maketitle

\section{Introduction and main results}

Let $(M,g)$ be a compact Riemannian manifold with boundary.  Let $\Phi^t$ be the geodesic flow on $TM$, where for each $(x,\xi)\in M$, $t\mapsto \Phi^t(x,\xi)$ is defined over its maximal interval containing $t=0$, in particular this interval is allowed to be the zero point only. Let $SM$ be the unit tangent bundle. Then $\partial SM$ represents all elements in $SM$ with a base point on $\bo$. 

Denote
\be{a1}
\partial_\pm  SM = \left\{(x,\xi)\in  \partial SM; \; \pm\langle \nu,\xi\rangle < 0\right\},
\ee
where $\nu$ is the unit  interior  normal, $\langle \cdot,\cdot\rangle$ and  stands for the inner product. The \textit{scattering relation} 
\be{sig}
\Sigma : \partial_-SM \to  \overline{\partial_+SM}
\ee
is defined by $\Sigma(x,\xi) = (y,\eta) = \Phi^\mathcal{L}(x,\xi)$, where $\mathcal{L}>0$ is the first moment, at which the (unit speed)  geodesic through $(x,\xi)$ hits $\bo$ again. If such an $\mathcal{L}$ does not exist, we formally set $\mathcal{L}=\infty$ and  we call the corresponding geodesic \textit{trapped}. This defines also $\mathcal{L}(x,\xi)$ as a function  $\mathcal{L}: \partial_-SM \to [0,\infty]$. Note that $\Sigma$ and $\mathcal{L}$ are not necessarily continuous. 

It is convenient to think of $\Sigma$ and $\mathcal{L}$ as defined on the whole $\partial SM$ with $\Sigma=\Id$ and $\mathcal{L}=0$ on $\overline{\partial_+SM}$. 

We parametrize the scattering relation in a way that makes it independent of pulling it back by diffeomorphisms fixing $\bo$ pointwise. Let $\kappa_\pm : \partial_\pm  SM \to B(\bo)$ be the orthogonal projection onto the (open) unit ball tangent bundle that extends continuously to the closure of $\partial_\pm  SM$. Then $\kappa_\pm$ are homeomorphisms, and we set 
\be{sig1}
\sigma = \kappa_+\circ \Sigma\circ \kappa^{-1}_- :\overline{B(\bo)} \longrightarrow \overline{B(\bo)}, \quad \ell = \mathcal{L}\circ \kappa^{-1}_-:  \overline{B(\bo)}\longrightarrow [0,\infty].
\ee
According to our convention, $\sigma=\Id$, $\ell=0$ on $\partial( \overline{B(\bo)} ) = S(\bo)$. We equip $\overline{B(\bo)}$ with the relative topo\-logy induced by $T(\bo)$, where neighborhoods of boundary points (those in $S(\bo)$) are given by half-neighborhoods. 

Let $\D$ be an open subset  of $\overline{B(\bo)}$. 
The \textit{lens rigidity} question we study in this paper is the following: 

\medskip
\textit{Given $M$ and $g|_{T(\bo)}$, do $\sigma$ and  $\ell$, restricted to $\D$, determine $g$ uniquely, up to a pull back of a diffeomorphism that is identity on $\bo$?}

\medskip

More generally, one can ask whether one can determine the topology of $M$ as well. 
One motivation for the lens rigidity problem is the
study of the inverse scattering problem for metric perturbations of the Laplacian. Suppose that
we are in Euclidean space equipped with a Riemannian metric which is Euclidean
outside a compact set. The inverse problem is 
to determine the Riemannian metric from 
the scattering operator, which is a Fourier integral operator, if
the metric is non-trapping (see \cite{Gu}). 
It was proven in \cite{Gu} that from the wave front set of
the scattering
operator, one can determine, under some conditions on the metric including non-trapping,
the scattering relation on the boundary of a
large ball. This uses high frequency information of the scattering
operator. In the semiclassical setting, Alexandrova has shown
that the scattering operator associated
to potential and metric perturbations of the Euclidean Laplacian
is a semiclassical Fourier integral operator that quantizes  the
scattering relation \cite{A1}, \cite{A2}. The scattering relation is also encoded in the hyperbolic Dirichlet to Neumann map on $\bo$. 
Lens rigidity is also
considered in \cite{PoR} in the study of the AdS/CFT duality and holography,
namely the idea that the ``bulk'' space-time can be captured by conformal field theory on a ``holographic screen''. The lens rigidity problem
appears also naturally when considering rigidity questions in Riemannian
geometry \cite{C1,C2}.

The lens rigidity problem is also closely related to the {\em boundary rigidity}
problem.  Denote by $\rho_g$ the  distance
function in the metric $g$. The boundary rigidity problem consists
of whether $\rho_g(x,y)$, known for all $x$, $y$ on $\partial M$, determines the metric uniquely. It is clear that any isometry which is the identity at the boundary will give rise to the same distance functions on the boundary. Therefore, the natural question is whether this is the only obstruction to unique identifiability of the metric.
The boundary distance  function only takes into account the  shortest paths and it is easy to find counterexamples  where $\rho_g$ does not carry any information about certain open subset of $M$, 
 so one needs to pose some restrictions on the metric. One such condition
is simplicity of the metric.

\begin{definition} We say that the Riemannian metric $g$ is {\em simple} in $M$, if $\partial M$ is strictly convex w.r.t.\ $g$, and for any $x\in  M$,  the exponential map $\exp_x :\exp^{-1}_x(M) \to M$ is a diffeomorphism.
\end{definition}

The manifold $(M,g)$ is called {\em boundary rigid} if one can determine the
metric  (and more generally, the topology) from the
boundary distance function up to an isometry which is the identity at the boundary.
It is a conjecture of Michel \cite{M} that the simple manifolds are boundary
rigid. This has been proved recently  in two dimensions \cite{PU}, for subdomains
of Euclidean space \cite{Gr} or for metrics close to Euclidean \cite{BI}, or symmetric spaces of negative curvature 
\cite{BCG}.   It was shown in \cite{SU-rig} that metrics   a~priori close to a  metric in a generic set, which includes real-analytic
metrics, are boundary rigid. For other local results see \cite{CDS}, \cite{E},
\cite{LSU}, \cite{SU1}. The lens rigidity problem is equivalent to
the boundary rigidity problem if the manifold is simple \cite{M}. 

Of course there is more information in the lens rigidity problem if the manifold is not simple. Even so, 
the answer to the lens rigidity problem, even when $\D ={B(\bo)}$, is negative,
as shown by the examples in \cite{CK}. Note that in these examples the
manifold is trapping, that is, there are geodesics of infinite length. The natural conjecture is that the scattering relation for non-trapping manifolds determines the metric uniquely up to the natural obstruction (\cite{U}). There are very few
results about this problem when the manifold is not simple.  Croke has shown that if a manifold
is lens rigid, a finite quotient of it is also lens rigid \cite{C2}.

This kind of data is  overdetermined. At least when $(M, g)$ is simple, knowledge of $x$, $y$, $\ell$ (the graph of the boundary distance function) determines $\xi$, $\eta$ uniquely \cite{M}. This can be extended, under some non-degenerate assumptions,  to the general case, see section~\ref{sec_jet}. Instead of the boundary distance function, we need to know the multiple travel times between boundary points. Also, for simple manifolds, $\sigma$ known on the whole $\partial_-SM$ determines $\ell$ uniquely. 

We  redefine $\sigma$ in a way that removes the need to know $g$ on $T(\bo)$. Denote by $T^0(\bo)$ the tangent bundle $T(\bo)$ considered as a conic set, i.e.,  vectors with the same direction in $T(\bo)$ are identified. 
For any metric $g|_{T(\bo)}$, $T^0(\bo)\setminus 0$ is isomorphic to the unit tangent bundle $S(\bo)$ (in the metric $g$) but has the advantage to be independent of the choice of $g$. 
Given $0\not=\xi'\in \overline{B_x(\bo)}$, we set
\be{l}
\lambda = |\xi'|_g\in [0,1], \quad \theta = \xi'/|\xi'|_g \in T_x^0(\bo),
\ee
i.e., $\lambda$ and $\theta$ are polar coordinates of $\xi'$. If $\xi'=0$, then $\theta$ is undefined. If $\xi' = \kappa_\pm(\xi)$, 
knowing $\lambda$, $\theta$ is equivalent to knowing the angle that $\xi$ makes with the boundary, and the direction of the tangential projection of the same vector. Given two metrics $g$ and $\hat g$ on $M$, and $(x,\xi')\in \overline{B(\bo)}$, $(x,\hat \xi')\in  \overline{\hat B(\bo)}$, where $\hat B(\bo)$ is related to $\hat g$, we say that $\xi'\equiv \hat \xi'$ iff $|\xi'|_g = |\hat\xi'|_{\hat g}$, and $\xi'=s\hat \xi'$ for some $s>0$. In other words, we require that $\xi'$ and $\hat\xi'$ have the same polar coordinates \r{l}. Note that this induces a homeomorphism $\overline{B(\bo)} \mapsto\overline{\hat B(\bo)}$ given by $\xi' \mapsto |\xi'|_g\xi'/|\xi'|_{\hat g}$ if $\xi'\not=0$, $0\mapsto 0$. 
With that identification of $B(\bo)$ for different metrics, it makes sense to study $\sigma$ restricted to the same set $\D$  for a family of metrics, and in particular, an a priori knowledge of $g$ on $T(\bo)$ is not needed to define $\D$ and $\sigma$ on it. If $\sigma(x,\xi')=(y,\eta')$, we just think of $\xi'$ and $\eta'$ as expressed in the polar coordinates \r{l}. Also, the notion of $\D$ being open is independent of $g$. 

A linearization of the boundary rigidity problem and the lens rigidity problem, see section~\ref{sec_lin}, 
is the following integral geometry problem. Given a family of geodesics $\Gamma$ with endpoints on $\bo$, we define the ray transform 
\be{ray}
I_\Gamma f(\gamma) = \int \langle f(\gamma(t)), \dot\gamma^2(t) \rangle \, \d t, \quad \gamma\in \Gamma,
\ee
of  symmetric 2-tensor fields $f$ (playing the role of the variation of the metric $g$), where   $\langle f,\theta^2\rangle$ is the action of $f$ on the vector $\theta$. Locally,  $\langle f,\theta^2\rangle = f_{ij}\theta^i\theta^j$. 
We will omit the subscript  $\Gamma$ to denote an integral over a chosen geodesic, or over all geodesics. Any such $f$ can be  decomposed orthogonally into a potential part $dv$ and a solenoidal one $f^s$ (see section~\ref{sec_4}), and $I$ vanishes on potential tensors. 
The linearized boundary rigidity or lens rigidity problem then is the following:  can we recover uniquely  the solenoidal projection $f^s$ of $f$ from its ray transform? If so, we call $I_\Gamma$ \textit{s-injective}. 

S-injectivity of $I$  was proved in \cite{PS} for metrics with negative curvature, in \cite{Sh, Sh-sib, D, Pe} for metrics with small curvature and in \cite{SU} for Riemannian surfaces with no focal points. A conditional and non-sharp stability estimate for metrics with small curvature is also established in \cite{Sh}. In
\cite{SU-Duke}, we proved stability estimates for s-injective simple metrics, see \r{401_1};   and sharp estimates about the recovery of a 1-form $f=f_j dx^j$ and a function $f$ from the associated $I f$. These  stability estimates  were used in \cite{SU-Duke} to prove local uniqueness for the boundary rigidity problem near any simple metric $g$ with s-injective $I$.  
In \cite{SU-rig}, we showed that the simple metrics $g$ for which $I$ is s-injective is generic, and applied this to the boundary rigidity problem. We note that in all the above mentioned results the metric has no conjugate points. On the other hand, in \cite{SU-AJM} we proved generic \text{s-injectivity} for a class of non-simple manifolds described below. Within that class, the boundary is not necessarily strictly convex, we might have conjugate points on the metric, the manifold might  be
trapping, and we have partial or incomplete information, i.e., we do not  know the scattering relation or the ray transform for all sets of geodesics.

Given $(x,\xi)\in \D$, let $\gamma_{\kappa_-^{-1} (x,\xi) }$ denote the geodesic issued from $\kappa_-^{-1} (x,\xi)$ with endpoint  $\pi(\sigma(x,\xi))$, where $\pi$ is the natural projection onto the base point. With some abuse of notation, we define
\[
I_\D(x,\xi) = I(\gamma_{ \kappa_-^{-1} (x,\xi)}), \quad (x,\xi)\in\D.
\]

\begin{definition}   \label{def_comp} 
We say that $\D$ is \textbf{complete} for the metric $g$, if for any $(z,\zeta)\in T^*M$ there exists a maximal in $M$, finite length unit speed geodesic $\gamma :[0,l] \to M$ through $z$, normal to $\zeta$, such that 
\begin{align}  \label{i}
&\left\{(\gamma(t), \dot\gamma(t)); \;0\le t\le l  \right\}   \cap S(\bo)  \subset  \D,\\ \label{ii}
&\text{there are no conjugate points on $\gamma$.}
\end{align}
\noindent We call the $C^k$ metric $g$ \textbf{regular}, if a complete set $\D$ exists, i.e., if $\overline{B(\bo)}$ is complete.
\end{definition}

If $z\in\bo$ and $\zeta$ is conormal to $\bo$, then $\gamma$ may reduce to one point. 
Since \r{i} includes  points where $\gamma$ is tangent to $\bo$, and $\sigma=\Id$, $\ell=0$ there, knowing $\sigma$ and $\ell$ on them provides no information about the metric $g$. On the other hand, we require below that $\D$ is open, so the purpose of \r{i} is to make sure that we know $\sigma$ near such tangent points.

\begin{definition} We say that $(M,g)$ satisfies the \textbf{Topological Condition (T)} if any path in $M$ connecting two boundary points  is homotopic to a polygon $c_1\cup \gamma_1 \cup c_2\cup\gamma_2\cup\dots \cup\gamma_k \cup c_{k+1}$  with the properties that for any $j$,

(i) $c_j$ is a path on $\bo$;

(ii) $\gamma_j :[0,l_j]\to M$ is a  geodesic lying  in $\Mint$ with the exception of its endpoints and is transversal to $\bo$ at both ends; moreover, $\kappa_-(\gamma_j(0), \dot\gamma_j(0)) \in  \D$; 

\end{definition}

Notice that (T) is an open condition w.r.t.\ $g$, i.e., it is preserved under small $C^2$ perturbations of $g$. 

We showed in \cite{SU-AJM} that if $\D$ is complete, then $I_\D f$ recovers the singularities of $f^s$. Next, see also Theorem~\ref{thm_2} below,  under the same conditions, and assuming (T) as well, $I_\D$ is s-injective for real-analytic metrics, and if $k\gg2$, also for generic metrics. 

To define the $C^k(M)$ norm in a unique way, and to make sense of real analytic $g$'s, we choose and fix a finite real analytic atlas on $M$.

Theorem~\ref{thm_1} below says, loosely speaking, that for the classes of manifolds and metrics we study, the uniqueness question for the non-linear lens rigidity problem can be answered locally by linearization. This is a non-trivial implicit function type of theorem however because our success heavily depends on the a priori stability estimate that the s-injectivity of $I_\D$ implies, and the latter is based on certain hypoelliptic properties of $I_\D$, as shown in \cite{SU-AJM}, see \r{401_1}. 
We work with two metrics $g$ and $\hat g$; and will denote objects related to $\hat g$ by $\hat \sigma$, $\hat \ell$, etc. Note that (T) is not assumed in the first theorem.

\begin{theorem}  \label{thm_1}
Let $g_0\in C^k(M)$ be a regular Riemannian metric on $M$ with $k\gg2$ depending on $\dim(M)$ only. Let $\D$ be open and complete for $g_0$, and assume that there exists $\D'\Subset\D$ so that $I_{g_0,\D'}$ is s-injective. Then there exists $\eps>0$, such that for any two metrics $g$, $\hat g$ satisfying  
\be{thm3}
\|g-g_0\|_{C^k(M)}+ \|\hat g-g_0\|_{C^k(M)}\le\eps,
\ee
the relations 
\[
\sigma = \hat\sigma, \quad \ell = \hat\ell \quad \text{on $\D$}
\]
imply that there is a $C^{k+1}$ diffeomorphism $\psi :M \to M$ fixing the boundary such that
\[
\hat g=\psi^*g.
\]
\end{theorem}

Next theorem is a version of \cite[Theorem~3]{SU-AJM}. It states that the requirement that $I_{g_0,\D'}$ is \mbox{s-injective} is a generic one for $g_0$. 

\begin{theorem}  \label{thm_2}
Let $\mathcal{G}\subset C^k(M)$, with $k\gg2$ depending on $\dim(M)$ only,  be an open set of regular Riemannian metrics on $M$ such that (T) is satisfied for each one of them. Let the  set $\D' \subset \overline{B(\bo)}$ be open and complete for each $g\in\mathcal{G}$. Then there exists  an open and dense 
subset  $\mathcal{G}_s$  of $\mathcal{G}$ such that $I_{g,\D'}$ is s-injective for any $g\in \mathcal{G}_s$.
\end{theorem}

Theorems~\ref{thm_1} and \ref{thm_2} combined imply that there is local uniqueness, up to isometry, near a generic set of regular metrics. 

\begin{corollary}  \label{cor_1}
Let $\D'$, $\mathcal{G}$, $\mathcal{G}_s$ be as in Theorem~\ref{thm_2}, and let $\D\Supset \D'$ be open and complete for any $g\in \mathcal{G}$. 
Then the conclusion of Theorem~\ref{thm_1} holds for any $g_0\in\mathcal{G}_s$. 
\end{corollary}

In section~\ref{sec_jet} we also prove that one can recover the jet of $g$ on $\bo$ in boundary normal coordinates from $\sigma$, $\ell$ under a non-conjugacy assumption. There is no generic assumption, and the recovery is actually explicit.

\begin{remark}  \label{rem_1}
Condition (T) in Theorem~\ref{thm_2}, and Corollary~\ref{cor_1} in some cases can be replaced by the assumption that $(M,g)$ can be extended to $(\tilde M,\tilde g)$ that satisfies (T). One such case is if $(\tilde M,\tilde g)$ is a simple manifold, and we study $\sigma$, $\ell$ on its maximal domain, i.e., $\D= \overline{B(\bo)}$. In particular, we get local generic lens rigidity for subdomains of simple manifolds when $\D$ is maximal. See section~\ref{sec_rem1} for more details.
\end{remark}

\section{Preliminaries}
We allow the geodesics to have segments on $\bo$, then they are called geodesics if they satisfy the geodesic equation in a half-neighborhood of each boundary point. Such segments are included in determining the maximal interval, where $\Phi^t(x,\xi)$ is defined. If $(\tilde M,\tilde g)$ is any extension of $(M,g)$ (a Riemannian manifold of the same dimension of which $M$ is a submanifold), then any geodesic in $\tilde M$ restricted to $M$ is a geodesic in $M$. For the maximal geodesics in $M$ we have the following property that indicates that the property of $\gamma$ to be maximal in $M$ does not change under such extensions. 

\begin{lemma}  \label{lemma-1}
Let $\gamma(t)$, $0\le t\le l$, $0\le l<\infty$, be a maximal geodesic in $M$. Let $(\tilde M, \tilde g)$ be any $C^{1,1}$ extension of $(M,g)$. Then there is no interval $I\supset [0,l]$ strictly larger than $[0,l]$, such that  $\gamma$ can be extended as a geodesic in $\tilde M$ for $t\in I$, and $\{\gamma(t); \; t\in I\}\subset M$. 
\end{lemma}

\begin{proof}
Suppose that there is such $I$. Without loss of generality, we may assume that $I\supset [0,l+\delta]$, $\delta>0$. Then   $\gamma$ solves the geodesic equation for $t\in [0,l+\delta]$, with the Christoffel symbols $\Gamma_{ij}^k$ depending on $g$ and its first derivatives restricted to $M$. Since $\tilde g$ is continuous and has continuous first derivatives across $\bo$, their restriction to $M$ depends on $g$ only. Therefore, $\gamma$ is a geodesic in $M$ for $t\in I$, and this is a contradiction.
\end{proof}

In \cite{SU-AJM} we studied geodesics originating from points outside $M$ given some extension $(\tilde M,\tilde g)$. We  connect the notion of $\D$ being open with the analysis in \cite{SU-AJM} by the following.

\begin{lemma}   \label{lemma_open}
Let $\tilde M$ be an extension of $M$, and let $\tilde g$ be a $C^{1,1}$ extension of $g$ on $\tilde{M}$. Let $\gamma_0 :[0,l]\mapsto \tilde{M}^\text{\rm int}$ be a unit speed geodesic  with endpoints in $\tilde{M}^\text{\rm int}\setminus M$ such that 
\be{a3}
\big\{(\gamma_0(t),\dot\gamma_0(t));\; t\in [0,l]\big\} \cap \overline{\partial_- SM} \subset  \kappa_-^{-1}(\D).
\ee
Then there exists a neighborhood $W$ of $(x_0,\xi_0)=(\gamma_0(0), \dot\gamma_0(0))$ such that any geodesic $\gamma$ with initial conditions in $W$ and the same interval of definition $t\in [0,l]$ still satisfies \eqref{a3}. 

Moreover, if $\tilde g$ belongs to a class of extensions satisfying $\|\tilde g\|_{C^{1,1}}\le A$ with some $A>0$, then $W$ can be chosen independently of $\tilde g$. 
\end{lemma}

\begin{proof}  
Let $E_0$ be the  l.h.s.\ of \r{a3}. Then $E_0$ is compact in $S\tilde M$.  
Given any neighborhood  $U$ of $E_0$ (in $S \tilde M$), the set $\partial SM\setminus U$ is compact, therefore there exists a neighborhood $W$ of $(x_0,\xi_0)$ such that the geodesic flow through $W$ for $0\le t\le l$ will miss that set, therefore, its only common points with $\overline{\partial_- SM}$ must be in $U$. To
construct $U$, we first choose a neighborhood $U_{x,\xi}$ in $S \tilde M$ of each point $(x,\xi)\in \kappa_-^{-1}(\D)$ so that $U_{x,\xi}\cap \overline{\partial_- SM}\subset  \kappa_-^{-1}(\D)$. This is easy to do in local coordinates because $\D$ is open. Then we set $U = \cup_{(x,\xi)\in \kappa_-^{-1}(\D)} U_{x,\xi}$.

To prove the second part, we use the theorem of continuous dependence of solutions of an ODE, over a fixed interval, on the initial conditions and on the coefficients of the ODE. As long the Lipschitz constant related to the generator of the geodesic flow  is uniformly bounded, we can choose $W$ uniformly w.r.t.\ $\tilde g$. 
\end{proof}

\begin{lemma}  \label{lemma_a3}
Let $(M,g)$, $(\tilde M,\tilde g)$, $\gamma_0$ be as in  Lemma~\ref{lemma_open}. Let $H$ be a hypersurface through $x_0=\gamma_0(0)$ transversal to $\gamma_0$, and set $\mathcal{H} = SM\cap \pi^{-1}(H)$.  Then there exists a small enough neighborhood $U$ of $(x_0,\xi_0)$ in $\mathcal{H}$, so that the geodesics issued from $U$ with interval of definition $t\in[0,l]$ satisfy \r{a3} and are transversal to $\bo$ at all common points with it except for a closed set of initial conditions  in $U$ of measure zero. 
\end{lemma}

\begin{proof}
We only need to prove the statement about the measure zero set. We study the geodesic flow on $SM$ issued from $U$. The corresponding geodesic is tangent to $\bo$ at some point $x$ if the corresponding integral curve in the phase space is tangent to $\partial SM$ at $\pi^{-1}(x)$. By Sard's theorem, this happens only on a closed set of measure zero. 
\end{proof}

Those lemmas will be used to reformulate the results in \cite{SU-AJM} in the situation in the paper. In \cite{SU-AJM}, we extended the geodesics slightly outside $M$ and parametrized them by initial points and directions on surfaces $H_m$ transversal to them, see section~\ref{sec_4.2}, instead by points and directions on $\bo$. This can be done, because  when studying $I_\Gamma$, the metric $g$ is known and can be extended outside $M$ in a known way. The reason for doing this was to  prevent working with a parametrization, where the geodesics can be tangent to the surface, which is the case with $\bo$. We can still do this even for the lens rigidity problem, using the boundary recovery result in section~\ref{sec_jet} below, see Proposition~\ref{pr_1}. We prefer however to parametrize the scattering relation by points on $\bo$ and corresponding directions. This does not preserve the smooth structure of the previous parametrization but the lemmas above show that it preserves the topology, at least.

\section{Recovery of the jet of $g$ on $\bo$.}  \label{sec_jet} 
If $\bo$ is convex in a neighborhood of some point $x_0$, then it is known that the jet of $g$ near $x_0$ in boundary normal coordinates is uniquely determined by studying the ``short'' geodesics connecting $x_0$ and $y$, and then letting $y\to x_0$, see \cite{LSU}. This argument does not provide explicit recovery however. If $\bo$ is strictly convex near $x_0$, then there is also a conditional Lipschitz stability estimate, see \cite{SU-rig}. Below we show that we can do this also without the convexity assumption by studying ``long geodesics'' that do not converge to a point. A non-conjugacy condition will be imposed. Note that there is no generic assumption in the theorem below, and the proof gives explicit recovery.
 
Let $(\tilde M, \tilde g)$ be a smooth extension of $(M,g)$. 
Let $x_0\in M$, $y_0\in M$ be endpoints of a   geodesic $\gamma_0: [0,1]\to\in M$, and assume that $x_0$, $y_0$, are not conjugate points on $\gamma_0$.  Then there exist neighborhoods $U\ni x_0$, $V\ni y_0$ in $\tilde M$, such that for $x\in U$, the exponential  map $\n(\dot\gamma_0(0))\ni \xi\mapsto y= \exp_x\xi\in V$ is a local diffeomorphism, therefore it has a smooth inverse $V\ni y\mapsto \xi :=\exp_x^{-1}y\in T_x\tilde M$  such that $\exp_{x_0}^{-1}y_0 =\dot\gamma_0(0)$. The geodesic $\gamma_{x,\xi}:[0,1]\to \tilde M$ then connects $x$ and $y$ and is unique among the geodesics issued from $x$ in directions close enough to $\xi$. One can also define a smooth travel time function $\tau(x,y)$ on $U\times V$ by $\tau(x,y) = |\exp_x^{-1}y|_g$ with $\tau(x_0,y_0)=|\dot\gamma(0)|_g$. If there are no conjugate points on $\gamma_0$, then $\tau$ locally minimizes the distance but in general, this is no longer true. On the other hand,  $\tau$ is always a critical value of the length functional and of the energy functional. 

In the  situation described in the paragraph above, we will call $y\in V$ \textit{visible} from $x\in U$, if $\gamma_{x,\xi}\subset M$. In the next theorem, given $(x,\xi)\in TM$, we call $y\in M$ \textit{reachable} from $(x,\xi)$, if there exists $s\ge0$, such that $\gamma_{x,\xi}(s)\in M$ for $t\in [0,s]$, and  $\gamma_{x,\xi}(s)=y$. 

The so defined function $\tau$ solves the eikonal equation $|\text{grad}_x\,\tau|_g^2=1$ in $U$ 
despite the possible existence of  pairs of conjugate points not in $U\times V$. Indeed, clearly, $d\tau(\gamma_{x,\xi}(t),y)/dt|_{t=0}=-|\xi|_g$. 
Next, by the Gauss lemma, the $x$-derivatives of $\tau(x,y)$ in directions perpendicular to $\xi$  vanish. Therefore, $\langle \xi,\text{grad}_x\, \tau\rangle = -|\xi|_g$, and the two vectors are parallel, therefore, $\text{grad}_x\, \tau=-\xi/|\xi|_g $, and $| \text{grad}_y\, \tau|_g^2=1$. 
Also, if $\eta = \text{grad}_y\tau$, then for some $l\ge0$, 
$\Phi^l(x,\xi) = (y,\eta)$. In particular, $\eta' = \text{grad}_y'\tau$, where the prime stands for a tangential projection.  If $\sigma$ is locally known, then integrating  the known $\eta'=\eta'(y)$ along a curve on $\bo$ connecting $y$ and $y_0$ recovers uniquely $\tau(x_0,y)$ up to a constant. 

Assume further that $\gamma_0$ is transversal to $\bo$ at both ends and does not touch $\bo$ elsewhere. Then 
\be{tau}
\sigma(x,-\text{grad}_x'\tau(x,y)) = (y,\text{grad}_y'(x,y)), \quad
\ell(x,-\text{grad}_x'\tau(x,y)) = \tau(x,y),
\ee
Therefore, at least in this non-degenerate situation, knowledge of $\tau(x,y)$ recovers uniquely $\sigma$, $\ell$ locally. 

The travel (arrival) times are widely used in the applied literature. Assume that $g$ on $T(\bo)$ is known and fixed. Fix $x$, $y$ on $\bo$. We call the number $\tau\ge0$  a {\em travel time} between $x$ and $y$, if there exists a geodesic of length $\tau$ connecting $x$ and $y$. Then we have a map that associates to any $(x,y)\in\bo\times\bo$ a subset of $(0,\infty]$. In  the situation above, $\tau(x,y)$ is one of the possible travel times between $x$ and $y$.

\begin{theorem}  \label{thm_jet}
Let $(M,g)$ be a compact Riemannian manifold with boundary. Let $(x_0,\xi_0)\in S(\partial M)$ be such that the maximal geodesic $\gamma_0$ through it is of finite length,  and assume that $x_0$ is not conjugate  to any point in $\gamma_0\cap \bo$.  If $\sigma$ and $\ell$ are known on some neighborhood of $(x_0,\xi_0)$, then  the jet of $g$ at $x_0$ in boundary normal coordinates is determined uniquely. 
\end{theorem}

\begin{proof}
To make the arguments below more transparent, assume that the geodesic $\gamma_0$ issued from $(x_0,\xi_0)$ hits $\bo$ for the first time transversally at $\gamma_0(l_0)=y_0$, $l_0>0$. If $l_0=0$, the results is known \cite{LSU}, as we pointed out above (but this proof applies to this case as well). 
 Then $y_0$ is the only point on $\bo$ reachable from $(x_0,\xi_0)$, and $x_0$, $y_0$ are not conjugate points on $\gamma_0$ by assumption. Assume also that $\bo$ is strictly convex at $x_0$, when viewed from the interior, i.e., the second fundamental form at $x_0$ is strictly negative. Then there is a half neighborhood $V$ of $x_0$ on $\bo$ visible from $y_0$ \cite{Sh-sib}. The latter is not always true for example when $\gamma_0$ is tangent to $\bo$ of infinite order at $x_0$. 

Choose local boundary normal coordinates near $x_0$ and $y_0$, and let $g_0$ be the Euclidean metric in each of them w.r.t.\ to the so chosen coordinates. In boundary normal coordinates, $\sigma$, $\ell$ determine uniquely $\Sigma$, $\mathcal{L}$. 
We can then consider a representation of $\Sigma$, denoted by $\Sigma^\sharp$ below,  defined locally on $\R^{n-1}\times S^{n-1}$, with values on another copy of the same space. If $(x,\theta)\in \R^{n-1}\times S^{n-1}$, then the associated vector at $x\in \bo$ is $\xi = \theta/|\theta|_g$; and $\Sigma^\sharp(x,\theta) =\Sigma(x,\xi)$. The same applies to the second component of $\Sigma^\sharp(x,\theta)$. Namely, if $(y,\eta) = \Sigma(x,\xi)$, then we set $\omega = \eta/|\eta|_{g_0}$, then $\Sigma^\sharp : (x,\theta) \mapsto (y,\omega)$. Similarly, we set $\mathcal{L}^\sharp(x,\theta) = \mathcal{L}(x,\xi)$. Let also $\theta_0$ and $\omega_0$ correspond to $\xi_0$ and $\eta_0$, respectively, where $\Sigma(x_0,\xi_0) = (y_0,\eta_0)$.

Set $\tau(x) := \tau(x,y_0)$, where $\tau$ is the travel time  function described above defined in a small neighborhood $U$ of $x_0$ in $\tilde M$. In the  normal boundary coordinates $x=(x',x^n)$ near $x_0$,  $g_{in}=\delta_{in}$, $\forall i$. Since $x_0$ and $y_0$ are not conjugate, for $\eta\in S_{y_0}M$ close enough to $\eta_0$,  the map $\eta\mapsto x\in\bo$ is a local diffeomorphism as long as the geodesic connecting $x$ and $y_0$ is not tangent to $\bo$ at $x$. Moreover, that map is known because it is determined by the inverse of $\Sigma$ near $(x_0,\xi_0)$. Similarly, the map $S^{n-1}\ni \omega \mapsto x$ is a local diffeomorphism and is also known. Then we know $(x,-\theta) = \Sigma^\sharp (y_0,-\omega)$, and we know $\mathcal{L}^\sharp(y_0,-\omega) = \mathcal{L}^\sharp(x,\theta) = \tau(x)$. 
Then we can recover $\text{grad}'\, \tau = -\theta'/|\theta|_g$.  Taking the limit $\omega\to\omega_0$, we recover $|\theta_0|^2_g = g_{\alpha\beta}\theta_0^\alpha\theta_0^\beta$. We use now the fact that a symmetric $n\times n$ tensor $f_{ij}$ can be recovered by knowledge of $f_{ij}v_k^iv_k^j$ for $N= n(n+1)/2$ ``generic'' vectors $v_k$, $k=1,\dots,N$; and such $N$ vectors exist in any open set on $S^{n-1}$, see e.g. \cite{SU-AJM}. Thus choosing appropriate  $n(n-1)/2$ perturbations of $\theta_0$'s, we recover $g(x_0)$. Thus, we recover $g$ in a neighborhood of $x_0$ as well; we can assume that $V$ is covered by that neighborhood.

Note  that we know now all tangential derivatives of $g$ in $V\ni x_0$. We showed above that $\tau$ solves the eikonal equation 
\be{b1}
g^{\alpha\beta} \tau_{x^\alpha} \tau_{x^\beta} +\tau_{x^n}^2 =1.
\ee

Next, in $V$, we know $\tau_{x^\alpha}$, $\alpha\le n-1$,  we know $g$, therefore by \eqref{b1}, we get $\tau_{x^n}^2$. It is easy to see that $\tau_{x^n}\le0$ on the visible part $V$, so we recover $\tau_{x^n}$ there. 
We therefore know the tangential derivatives of $\tau_{x^n}$ on $\bo$ near $x_0$.

Differentiate \eqref{b1} w.r.t.\ $x^n$ at $x=x_0$ to get
\be{b2}
\left[
\frac{\partial g^{\alpha\beta}}{\partial  x^n} \tau_{x^\alpha} \tau_{x^\beta} + 2g^{\alpha\beta}   \tau_{x^\alpha x^n} \tau_{x^\beta} + 2  \tau_{x^n x^n}   \tau_{x^n}\right]\bigg|_{x=x_0} =0.
\ee
Since $\gamma_0$ is tangent to $\bo$ at $x_0$, we have $\tau_{x^n}(x_0)=0$ by \r{b1}. The third term in the l.h.s.\ of \eqref{b2} therefore vanishes.  Therefore the only unknown term in \eqref{b2} is $G^{\alpha\beta} := \partial  g^{\alpha\beta}/{\partial  x^n}$ at $x=x_0$. Since $\tau_{x^\alpha}(x_0) =  - \xi_0$, using the fact that $\text{grad}\,\tau(x_0) = -\xi_0$ again,  we get that we have to determine $G^{\alpha\beta}$ from $G_{\alpha\beta}\xi_0^\alpha\xi_0^\beta$. This is possible if we repeat the construction and replace $\xi_0$ by a finite number of vectors near $\xi_0$, as above. So we get an explicit recovery of $\partial g/\partial x^n|_{\bo}$ in fact. 

Next, for $x\in V$ but not on $\partial V$,  we can recover $\tau_{x^nx^n}(x)$ by \r{b2} because $\tau_{x^n}(x)<0$. By continuity, we recover $\tau_{x^nx^n}(x_0)$, therefore we know $\tau_{x^nx^n}$ near $x_0$, and all tangential derivatives of the latter. 

We differentiate \r{b2} w.r.t.\ $x^n$ again, and as above, recover $\partial^2g/\partial (x^n)^2|_{\bo}$ near $x_0$. Then we recover $\partial^3\tau/\partial (x^n)^3$, etc.

Let us return to the general case. We will show first that the assumptions about $(x_0,\xi_0)$ are preserved under a small perturbation of that point. 
Let $U$ be a neighborhood of $(x_0,\xi_0)$ in $TM$. Since $\gamma_0\cap \bo$ consists of points that are not conjugate to $x_0$, there is an open set $W\supset\gamma_0\cap \bo$ in $\tilde M$ that stays away from the points conjugate to $x$ along the geodesics  $t\mapsto \exp_{x}t\xi$, $0\le t\le l_0$ if $(x,\xi)\in U$ and if $U$ is small enough. Here $l_0$ is the length of $\gamma_0$. 
On the other hand, the possible common points of those geodesics must be in $W$, if $U$ is small enough, as in the proof of Lemma~\ref{lemma_open}. 
By Lemma~\ref{lemma-1}, there exists $s>1$ close enough to $1$ so that $(x_0,s\xi_0)\in  U$, and  $\exp_{x_0}(s\xi_0)\not \in M$. The geodesics issued from $(x,\xi)\in U$ close enough to $(x_0,s\xi_0)$, and the same time interval $[0,l_0]$ of definition, still have endpoints outside $M$, therefore they are longer than the maximal segment in $M$. On the other hand, we showed that they meet $\bo$ at points that are not conjugate to $x_0$. We can now replace $s\xi_0$ by $\xi_0$ and rescale the corresponding geodesic. This shows that the assumption in the theorem is preserved under a small enough perturbation of $(x_0,\xi_0)$.

Let $\xi_\eps= \xi_0+\eps\nu$, $0<\eps\ll1$, where $\nu$ is the interior unit normal at $x_0$. In the coordinate system above, $\nu=(0,\dots,0,1)$. Let $\gamma_\eps$ be the geodesic issued from $(x_0,\xi_\eps)$ until it hits $\bo$ for the first time.   
By a compactness argument, we can choose a sequence $\eps_j \searrow 0$ such that $y_{\eps_j} \to y_0^*$ that can be different from $y_0$ but we still have that $y_0^*\in\bo$ and $y_0^*$ is reachable from $(x_0,\xi_0)$. Then 
$l_{\eps_j}\to l_0^*$ with some $l_0^*\ge0$. 
By assumption, $x_0$ and $y_0^*$ are not conjugate points on $\gamma_0$, and the function $\tau(x,y)$ is then  well-defined  near $(x_0,y_0^*)$ satisfying the eikonal equation. 
Next,  $\gamma_{\eps_j}$  connects $x_0$ and $y_{\eps_j}$, and only the endpoints are not in $\Mint$. The advantage now is that $\gamma_{\eps_j}$ hits $\bo$ at $x_0$ transversely but it might be tangent to $\bo$ at $y_{\eps_j}$. 

For a fixed $j$, by \cite{Sh-sib} (see  the proof of Lemma~2.3 there)  if $U$ is a small enough neighborhood of $x_0$ on $\bo$, then $U$ can be expressed as the disjoint union $U^+\cup H\cup U^-$, where $H$ is a hypersurface on $\bo$ through $x_0$, and at least one of the half-neighborhoods $U^\pm$ is visible from $y_{\eps_j}$, let us say that this is $U^+$. If $\gamma_{\eps_j}$ is transversal to $\bo$ at $y_{\eps_j}$, then even better, the whole $U$ is visible from $y_{\eps_j}$ if $U$ is small enough. If $n=2$, then $H$ reduces to the point $x_0$ and the modifications are obvious. We set $\tau(x) = \tau(x,y_{\eps_j})$. 
By the arguments in the special case above, we know $\tau(x)$ on $\bo$.

We first recover $g$ on $\bo$ in boundary normal coordinates. Let $\theta_j = \xi_{\eps_j}/|\xi_{\eps_j}|_{g_0}$, and $\theta_0 = \xi_{0}/|\xi_{0}|_{g_0}$. As in the beginning of the proof, we take a one-sided derivative at $x_0$, from $U^+$,  to recover $-\theta_{\eps_j} /|\theta_{\eps_j}|_g = -\theta'_0/|\theta_0|_g + O(\eps_j)$. Take the limit $\eps_j\to0$ to recover $|\theta_0|_g$. Thus we get $g$ at $x_0$ and therefore, near $x_0$.

As above, from the eikonal equation, satisfied in $U^+$, where $\tau$ is defined and known up to a constant, if $U^+$ is small enough, we recover $\tau_{x^n}^2$ on $U^+$, and therefore $\tau^2_{x^n}(x_0)$ by continuity. Since $\tau_{x^n}(x_0)= -\eps_j/|\xi_{\eps_j}|_g$, and therefore, $\tau_{x^n}<0$ near $x_0$, we recover $\tau_{x^n}(x_0)$, and therefore, the normal derivative of $\tau$ on $\bo$ on $U^+\cup H$, and the tangential derivatives of the latter. 

By \eqref{b2}, 
\[
G_{\alpha\beta}(x_0)\xi_0^\alpha\xi_0^\beta /|\xi_{\eps_j}|^2_g = \left(-2g^{\alpha\beta} \tau_{x^\alpha x^n} \tau_{x^\beta} +2\tau_{x^nx^n} \eps_j/|\xi_{\eps_j}|_g\right) \Big|_{x=x_0},
\]
where $G^{\alpha\beta} = \partial^{\alpha\beta}/\partial x^n$.  Since $\tau$ is smooth near $(x_0, y_0^*)$, then $\tau_{x^nx^n}(x_0)= \tau_{x^nx^n}(x_0,y_{\eps_j})$ remains bounded and even has a limit, as $j\to\infty$. Similarly the other terms above have a limit. 
Therefore, the second term on the r.h.s.\ above tends to zero, as $j\to\infty$,  and we recover  $G_{\alpha\beta}(x_0)\xi_0^\alpha \xi_0^\beta$.

To recover $\partial g/\partial x^n$ at $x_0$, we need to perturb $\xi_0$. 
As we showed above, $(x_0,\xi)$ satisfies the assumptions of the theorem, if $\xi$ is close enough to $\xi_0$.  We can therefore recover $G_{\alpha\beta}(x_0)\xi^\alpha\xi^\beta$ for $\xi$ in a neighborhood of $\xi_0$, which recovers $G_{\alpha\beta}(x_0)$, and therefore $\partial g/\partial x^n$ at $x=x_0$. 

To recover the higher order derivatives, we proceed in the same way. To recover $\partial^kg/\partial(x^n)^k$, we differentiate \r{b1} $k$ times, and solve for $\partial^kg^{\alpha\beta}/\partial(x^n)^k$ at $x=x_0$. The only unknown term in the r.h.s.\ will be $\partial^{k+1}\tau/\partial(x^n)^{k+1}$ but it will be multiplied by $\tau_{x^n}$ that equals $\eps_j/|\xi_{\eps_j}|_g$ at $x_0$. Then taking the limit $j\to\infty$ will recover $\partial^kg/\partial(x^n)^k$ as above. Then we recover $\partial^{k+1}\tau / \partial (x^n)^{k+1}$ needed for the next step, etc.
\end{proof}

\begin{remark} If $g$ has a finite smoothness $g\in C^k(M)$, then the proof above implies that we can recover (explicitly) $\partial^\alpha g|_{\bo}$ for $|\alpha|\le k-2$ in boundary normal coordinates.
\end{remark}

\section{Local interior rigidity; Proof of Theorem~\ref{thm_1}}  \label{sec_4} 
Given a symmetric  2-tensor $f= f_{ij}$,  the {\em divergence}  of $f$  is an 1-tensor $\delta f$  defined by  
$$
[\delta f]_i = g^{jk} \nabla_k f_{ij}
$$ 
in any local coordinates, where $\nabla_k$ are the covariant derivatives of the tensor $f$. Given an 1-tensor (a vector field or an 1-form that we identify through the metric) $v$, we denote by $dv$ the 2-tensor called symmetric differential of $v$:
$$
[d v]_{ij} = \frac12\left(\nabla_iv_j+ \nabla_jv_i  \right).
$$
Operators $d$ and $-\delta$ are formally adjoint to each other 
in $L^2(M)$.  
It is easy to see that for each smooth $v$ with $v=0$ on $\partial M$, we have $I(d v)(\gamma)=0$ for any geodesic $\gamma$ with endpoints on $\bo$. This follows from the identity  
\be{v}
\frac{\d }{\d  t}  \langle v(\gamma(t)), \dot\gamma(t) 
\rangle \allowbreak  = \allowbreak \langle
dv(\gamma(t)),  \dot\gamma^2(t) \rangle.
\ee 

It is known (see \cite{Sh} and  \r{9} below) that for $g$ smooth enough, each symmetric tensor $f\in L^2(M)$ admits unique orthogonal decomposition $f=f^s+d v$ into a {\em solenoidal}\/ tensor  $\mathcal{S}f :=f^s $ and a {\em potential}\/ tensor $\mathcal{P}f :=d v$, such that both terms are in $L^2(M)$, $f^s$ is solenoidal, i.e., $\delta f^s=0$ in $M$, and $v\in H^1_0(M)$ (i.e., $v=0$ on $\partial M$). In order to construct this decomposition, introduce the operator $\upDelta^s = \delta d$ acting on vector fields. This operator is elliptic in $M$, the Dirichlet problem satisfies the Lopatinskii condition, and has a trivial kernel and cokernel.  Denote by $\upDelta^s_D$ the Dirichlet realization of $\upDelta^s$ in $M$. Then
\begin{equation}  \label{9}
v = \left(\upDelta^s_D\right)^{-1}\delta f, \quad
f^s = f - d \left(\upDelta^s_D\right)^{-1}\delta f.
\end{equation} 
Therefore, we have
$$
\mathcal{P} = d \left(\upDelta^s_D\right)^{-1}\delta, \quad \mathcal{S} = \Id-\mathcal{P},
$$
and for any $g \in C^1(M)$, the maps
\be{10}
(\upDelta^s_D)^{-1}: H^{-1}(M) \to H_0^{1}(M),  \quad 
\mathcal{P},  \mathcal{S} :  L^2(M) \longrightarrow L^2(M)
\ee
are bounded and depend continuously on $g$, see \cite[Lemma~1]{SU-rig} that easily generalizes for manifolds. 
This  admits the following easy generalization:  for $s=0,1,\dots$, the resolvent above also continuously maps $H^{s-1}$ into $H^{s+1} \cap H_0^1$, similarly, $\mathcal{P}$ and $\mathcal{S}$ are bounded in $H^{s}$, if $g\in C^k$, $k\gg2$ (depending on $s$).  Moreover those operators depend continuously on $g$. 

Notice that even when $f$ is smooth and $f=0$ on $\partial M$, then $f^s$ does not need to vanish on $\partial M$. In particular, $f^s$, extended as $0$ to $\tilde M$, may not be solenoidal anymore. 
To stress on the dependence on the manifold, when needed, we will use the notation $v_M$ and $f^s_M$ as well. 

Operators $\mathcal{S}$ and $\mathcal{P}$ are orthogonal projectors. 
The problem about the s-injectivity of $I$, restricted to a subset of geodesics, then can be posed as follows: if $If=0$ on those geodesics, show that $f^s=0$, in other words, show injectivity on the subspace $\mathcal{S}L^2$ of  solenoidal tensors. 

We start with the proof of Theorem~\ref{thm_1}. We will split the proof into several steps.

\subsection{Choosing a suitable metric isometric to  $\hat g$.} \label{sec_4.1} 
Any two metrics such that one of them is a pull-back of the other under a diffeomorphism fixing the boundary pointwise, will be called below isometric. Such a diffeomorphism is necessarily $C^{k+1}$ if the metrics are $C^k$, see e.g., \cite{SU-rig}, and the norm of its derivatives are controlled by those of the two metrics, see \cite[Lemma~6]{SU-rig}.

We  first find a metric isometric to $\hat g$, that we denote by $\hat g$ again, so that the boundary normal coordinates related to $g$ and $\hat g$ coincide in some neighborhood of the boundary.

By \cite[Theorem~2.1]{CDS}, if $\eps\ll1$, there is $h\in C^{k-1}$  isometric to $\hat g$ so that $h$ is solenoidal w.r.t.\ $g$. Moreover,
\be{I2b}
\|h-\hat g\|_{C^{k-1}} \le C\eps,
\ee
and there is no need to replace $k$ by $k-1$ if we work in the $C^{k,\alpha}$ spaces.  By a standard argument, by a diffeomorphism that identifies normal coordinates near $\bo$ for $h$  and $g$, and is identity away from some neighborhood of the boundary, we find a third $\hat g_1$ isometric to $h$ (and therefore to $\hat g$), so that $\hat g_1=\hat g$ near $\bo$, and $\hat g_1=h$ away from some neighborhood of $\bo$ (and there is a region that $\hat g_1$ is neither). Then $\hat g_1-h$ is as small as $g-h$, more precisely,
\be{I2a1}
\|\hat g_1-h\|_{C^{l-3}} \le C\|g-h\|_{C^{l-1}}, \quad l\le k.
\ee
This follows from the fact that $\hat g_1 = \phi^*h$, with a diffeomorphism $\phi$ that is identity on the boundary, and
\be{phi}
\|\phi-\Id\|_{C^{l-2}}\le C\|\hat g-h\|_{C^{l-1}}.
\ee
Set
\be{ft}
f = h-g, \quad \tilde f = \hat g_1-g.
\ee
We aim to show that $f=\tilde f=0$. 
Estimate \eqref{I2a1}  implies
\be{I2a2}
\|\tilde f-f\|_{C^{l-3}} \le C\|f\|_{C^{l-1}}, \quad \forall l\le k.
\ee
By \r{thm3}, \eqref{I2b} and \eqref{I2a2},
\be{I2c}
\|f\|_{C^{k-1}} \le C\eps, \quad \|\tilde f\|_{C^{k-3}} \le C\eps.
\ee
By Theorem~\ref{thm_jet}, and the remark after it,
\be{I3}
\partial^\alpha \tilde f=0 \quad \text{on $\bo$ for $|\alpha|\le  k-5$}.
\ee

We have now two isometric copies of $\hat g$: the first one is $h$ that has the advantage of being solenoidal w.r.t.\ $g$; and the second one $\hat g_1$ that has the same jet as $g$ on $\bo$. We need both properties below to show that $g=h$, i.e., $f=0$ (or $g=\hat g_1$, i.e., $\tilde f=0$) but so far we cannot prove that $h=\hat g_1$. The next proposition shows that  $h$ and $\hat g_1$ are equal up to $O(\|f\|^2) = O(\|\tilde f\|^2)$.  

\begin{proposition}   \label{pr_dv}
Let $\hat g$ and $g$ be in $C^k$, $k\ge2$ and isometric, i.e., 
\[
\hat g=\psi^*g
\]
for some diffeomorphism $\psi$ fixing $\partial M$ pointwise. Set $f = \hat g-g$. Then there exists $v$ vanishing on $\partial M$, so that 
\[
f = 2dv +f_2,
\]
and for  $g$ belonging to any bounded set $U$ in $C^k$, there exists $C(U)>0$, such that
\[
\|f_2\|_{C^{k-2}}\le C(U)\|\psi-\text{\rm Id}\|^2_{C^{k-1}},
\quad 
\|v\|_{C^{k-1}} \le C(U)\|\psi-\text{\rm Id}\|_{C^{k-1}}.
\]
\end{proposition}

\begin{proof}
Extend $g$ to $\tilde M$ in such a way that the $C^k$ norm of the extension is bounded by $C\|g\|_{C^k(M)}$. 
Set $v(x)= \exp^{-1}_x(\psi(x))$ that is a well defined vector field in $C^1( M)$  if $\psi$ is close enough to identity in $C^1$ (it is enough to prove the theorem in this case only), and $v=0$ on $\partial $M. Set $\psi_\tau(x)  = \exp_x(\tau v(x))$, $0\le \tau\le1$. Let $g^\tau = \psi_\tau^* g$. Under the smallness condition above, $v$ is small enough in $C^1$, and therefore $\psi_\tau$ is close enough to identity in the $C^1(\tilde M)$ norm. Therefore, $\psi_\tau :M \to \psi_\tau(M)\subset \tilde M$ is a diffeomorphism. Next, $\psi_\tau$ fixes $\bo$ pointwise, therefore, $ \psi_\tau(M) =M$. 

The Taylor formula implies
\[
\hat g = g+ \frac{d}{d\tau}\Big|_{\tau=0} g^\tau +h = g+2d v+h,
\]
where
\[
|h|\le \frac12 \max_{\tau\in[0,1]} \Big|\frac{d^2 g^\tau}{d\tau^2}\Big|,
\]
and $2d v$ is the linearization of $g^\tau$ at $\tau=0$, see \cite{Sh}. To estimate $h$, write
\[
g^\tau_{ij} = g_{kl}\circ \psi_\tau \frac{\partial \psi^k_\tau}{\partial x^i} \frac{\partial \psi^l_\tau}{\partial x^j},
\]
and differentiate twice w.r.t.\ $\tau$. Notice that 
\[
\Big|\frac{\partial^2 \psi_\tau}{\partial\tau^2}\Big| \le C\|v\|^2_{L^\infty}, \quad
\Big|\frac{\partial^2 \nabla \psi_\tau}{\partial\tau^2}\Big| \le C\|v\|^2_{C^1}.
\]
This yields the stated estimate  for  $f_2$ for $k=2$. The estimates for $k>2$ go along similar lines by expressing the remainder $h$ in its Lagrange form, and estimating the  derivatives of $h$.
\end{proof}

We apply Proposition~\ref{pr_dv} to $h$ and $\hat g_1$ to get by \r{phi},
\be{I4}
\tilde f = f+2dv +f_2, \quad \|f_2\|_{C^{l-3}} \le  C\|f\|_{C^{l-1}}^2, \quad \forall l\le k.
\ee
In other words, $\tilde f^s=f$ up to $O(\|f\|^2)$. 

Next, with $g$ extended as above, we extend $\hat g_1$ so that $\hat g_1=g$ outside $M$. Then $g\in C^k$ and $\hat g_1\in C^{k-5}$ by \r{I3}.

\subsection{Reparametrizing the scattering relation.} \label{sec_4.2} 
We proceed with some preliminary work that would allow us to apply \cite[Theorem~2]{SU-AJM}. Assume first that the underlying metric is fixed to be $g_0$. Let $(\tilde M, \tilde g_0)$ be a $C^k$ extension as above. In \cite{SU-AJM}, the geodesics are extended to $\tilde M\setminus M$, parametrized by initial points, and corresponding directions on a finite collection $\{H_m\}$ of smooth connected hypersurfaces in $\tilde M$, having additional properties as explained below. Given two complete $\D'\Subset \D$, we will  construct such a family issued from a set $\D''$ with   $\D'\Subset\D'' \Subset\D$, that is also complete.

For any $(z_0,\zeta_0)\in T^*M$, including the case where $z_0\in\bo$, there is a maximal geodesic $\gamma_0$ through $z_0$ normal to $\zeta_0$, satisfying the conditions of Definition~\ref{def_comp}.  Let us assume that $\gamma_0$ is parametrized by $t\in [l^-,l^+]$, $\pm l^\pm\ge0$, and $\gamma(0)=z$, $\dot\gamma_0(0)=\xi_0$, with $\zeta_0$ conormal to $\xi_0$.  Since $\gamma_0$ is maximal in $M$, by Lemma~\ref{lemma-1}, for any $\delta_1 >0$ there exists $\delta \in (0,\delta_1)$ so that the extension $\tilde \gamma_0$ of $\gamma$ to $\tilde M$ corresponding to $t\in [l^- -\delta,l^+ +\delta]$ is well defined and has endpoints in $\tilde{M}^\text{\rm int} \setminus M$. By \r{i}, 
\[
\left\{(\gamma_0(t),\dot\gamma_0(t));\, l^-\le t\le l^+\right\} \cap S(\bo) \subset\D.
\]
The extension $(\tilde \gamma_0,\dot{\tilde{\gamma}}_0)$ may have additional point on $\bo$ corresponding to $l^- -\delta\le t\le l^-$, and $l^+\le t\le l^+ +\delta$. However, we choose $\delta_1\ll1$ so that they still belong to $\D$, and this is possible to do because $\D$ is open. Now, by Lemma~\ref{lemma_open}, any geodesic that is obtained by a small enough perturbation $(z,\xi)$ of the initial conditions $(z_0,\xi_0)$ of $\gamma_0$ at $t=0$, with the same interval $t\in [l^- -\delta,l^++\delta]$, will satisfy condition \r{a3}. Condition \r{ii} will also be satisfied by a perturbation argument, if $\delta$ is small enough as well. Now we can perturb $\tilde g_0$ in the $C^2(\tilde M)$ topology to ensure the same property. Note that if $\|g-g_0\|_{C^2(M)}\le \eps$, one can choose an  extension $\tilde g$ of $g$ to $\tilde M$ so that $\|\tilde g-\tilde g_0\|_{C^2(M)}\le C\eps$ with $\tilde g_0$ a fixed extension of $g_0$ as above. 

To summarize, we proved the following.

\begin{lemma}  \label{lemma_4}
Under the conditions of Theorem~\ref{thm_1}, for any $(z_0,\zeta_0)\in S^*M$, there exists a geodesic in the metric $\tilde g_0$ through $z_0$ normal to $\zeta_0$ with endpoints in $\tilde{M}^\text{\rm int} \setminus M$ so that conditions  \r{ii}, \r{a3} are satisfied. 

Moreover, if that geodesic has initial conditions $(z_0,\xi_0)$  at $t=t_0$ for some $t_0\in[0,l]$, and an interval of definition $0\le t\le l$, properties  \r{ii} and \r{a3}  remain true under small enough perturbations of $(z_0,\xi_0)$, and $\tilde g_0$ in $C^2(\tilde M)$.
\end{lemma}

Let us assume now that the underlying metric is $g$ as in the theorem, with $\eps\ll1$. Since $\overline{\D'}$ is compact, there are finitely many geodesics 
\[
\left\{\gamma_m(t); \; l_m^--\delta_m\le t\le l_m^+ +\delta_m\right\},
\]
with the following properties. If $\dot\gamma_m(0)=\xi_m\in S_{z_m}M$, then for any $m$  there exists  neighborhoods $U'_m\Subset U_m$ of $(z_m,\xi_m)$ in $SM$, such that  if $\Gamma_m$, $\Gamma'_m$ is the set of geodesic with initial conditions in $U_m$, respectively $ U_m'$, and the same interval of definition as $\gamma_m$, then for any $(x,\zeta)\in T^*M$ there is a geodesic $\gamma\in \cup \Gamma_m'$ so that  \r{ii}, \r{a3} are satisfied with $\D$ replaced by $\D'$; and all geodesics in $\cup\Gamma_m$ satisfy  \r{ii}, \r{a3} as well. Moreover, 
\be{Ga}
\cup \Gamma_m' \supset \kappa_-^{-1}(\overline{D'}),
\ee
where $\Gamma'_m$ is regarded as a point set in the phase space  $S\tilde M$ consisting of the points on all  integral curves.

We will parametrize  $\Gamma_m$, $\Gamma_m'$, with initial points outside $M$. Choose a family of finitely many small enough smooth hypersurfaces $\{H_m\}$ in $\tilde{M}^\text{\rm int} \setminus M$, each one transversal to $\gamma_m$. Without loss of generality, we can assume that all geodesics in $\Gamma_m$ can be extended in $\tilde{M}^\text{\rm int} \setminus M$ so that they intersect $H_m$ once and are transversal to $H_m$ as well. We still denote the set of the extended geodesics by $\Gamma_m$ and $\Gamma_m'$, respectively. 
Let $\mathcal{H}_m$ be the open subset of $\{ (x,\xi)\in SM; \, x\in H_m, \, \xi\not\in T_xH_m  \}$ formed by those $(x,\xi)$ that coincide with the left endpoint and the corresponding direction of some geodesic in $\Gamma_m$. Then their endpoints belong to $\tilde{M}^\text{\rm int} \setminus M$ again, and their length is a smooth function $l_m(x,\xi)>0$ (actually, we may even assume that $l_m$ is constant, if $U_m$ is small enough). 

We have 
\be{401}
\Gamma_m = \Gamma(\mathcal{H}_m) = \left\{\gamma_{x,\xi}(t); \; 0\le t\le l_m(x,\xi), \, (x,\xi)\in \mathcal{H}_m\right\},
\ee
where $\gamma_{z,\xi}$, as usual, is the geodesic with initial conditions $(x,\xi)$. 
We define $\mathcal{H}'_m$ in a similar way, related to $U'_m$.   
We consider also the geodesics in the metric $\hat g$ defined as in \r{401}. Then \r{ii}, \r{a3} hold for those geodesics, too, provided that $\eps\ll1$.

Combining the arguments above with Lemma~\ref{lemma_a3}, we have the following.

\begin{proposition}  \label{pr_1}
Let $g$ and $\hat g_1$ be as in the end of Section~\ref{sec_4.1}. Then, if $\eps\ll1$, 
\[
\Phi^{l_m(x,\xi)}(x,\xi) = \hat\Phi^{l_m(x,\xi)}(x,\xi), \quad \forall (x,\xi)\in \mathcal{H}_m, \; \forall m,
\]
where $\hat \Phi$ is related to $\hat g_1$. 
Moreover, $\cup \Gamma_m$ satisfies \r{ii}, \r{a3}, and $\cup \Gamma_m'$ is complete in the sense that \linebreak $N^*(\cup \Gamma_m')\supset T^*M$, and satisfies \r{Ga}; similarly   $\cup \hat\Gamma_m$ and $\cup \hat\Gamma_m'$ have the same properties. 
\end{proposition}

In other words, informally speaking, we pushed the boundary, where the scattering relation is defined, to a collection of hypersurfaces outside $M$, so that the corresponding geodesics are always transversal to them, and the endpoints are away from $\bo$.  

\begin{proof}[Proof of Proposition~\ref{pr_1}]
The proof is straightforward, if the geodesic issued from $(x,\xi)$, for $0\le l\le l_m(x,\xi)$, always intersects $\bo$ transversally. Observe that $\kappa_\pm=\hat \kappa_\pm$ because $g$ and $\hat g$ have the same normals on $\bo$, therefore $\sigma=\hat\sigma$ implies $\Sigma=\hat\Sigma$. 
The points $(x,\xi)\in \mathcal{H}_m$ where this transversality does not hold is a closed set of measure zero by  Lemma~\ref{lemma_a3}, and for such points, one can approximate with points outside this set, and to use the continuity of $\Phi^t$. 
\end{proof}

\subsection{Linearization near $g$.}   \label{sec_lin}
Now we are in the situation of \cite{SU-AJM}, see Theorem~2 there. Choose smooth functions $\alpha_m$ supported in $\mathcal{H}_m$ and equal to $1$ on $\mathcal{H}'_m$. Set $\alpha = \{\alpha_m\}$ and $I_{\alpha_m} = \alpha_m I$, more precisely,
\[
I_{\alpha_m} f(z,\xi) = \alpha_m(x,\xi)\int_{0}^{l_m(x,\xi)} \langle f(\gamma_{z,\xi}),\dot\gamma_{z,\xi}^2 \rangle \,\d t, \quad (z,\xi) \in \mathcal{H}_m.
\]
Also set
\[
I_\alpha = \{  I_{\alpha_m}\}, \quad N_{\alpha_m} = I_{\alpha_m}^* I_{\alpha_m}, \quad N_\alpha = \sum N_{\alpha_m},
\]
where the adjoint is taken w.r.t.\ the measure $|\langle \nu, \xi\rangle|\d\Sigma_{2n-2}$, where $\d\Sigma_{2n-2}$ is the induced measure on $\partial SM$ by the volume form, and $\nu$ is a unit normal to $H_m$.  

In \cite[Theorem~2]{SU-AJM} we showed that if for a fixed  $g$, $I_\alpha$ is s-injective, then we have the following a priori estimate
\be{401_1}
\|f^s\|_{L^2(M)} \le C\|N_\alpha f\|_{\tilde{H}^2(\tilde M)},
\ee
with a suitable norm $\|\cdot\|_{\tilde{H}^2(\tilde M)}$ so that $H^2(\tilde M)\subset \tilde{H}^2(\tilde M) \subset H^1(\tilde M)$. Moreover, \r{401_1} remains true under small $C^k$ perturbations of $g$ with a constant $C$ that can be choose uniformly. Note that (T) is not needed in \cite[Theorem~2]{SU-AJM}.

Fix $m$ and $(x,\xi)\in\mathcal{H}_m$. Let $\{\gamma_{x,\xi}(t)$, $0\le t\le l_m(x,\xi)\}$, and $\{\hat\gamma_{x,\xi}(t)$, $0\le t\le  l_m(x,\xi)\}$ be the geodesic issued from $(x,\xi)$ related to $g$ and $\hat g_1$, respectively. By Proposition~\ref{pr_1}, their endpoints and directions coincide. We reparametrize $\gamma_{x,\xi}$, $\hat \gamma_{x,\xi}$ so that $t\in[0,1]$; then they  have the same speeds $|\dot\gamma_{x,\xi}(t)|=  |\dot{\hat\gamma}_{x,\xi}(t)|= l_m(x,\xi)$ . 

Define the following variation of $\gamma_{x,\xi}$, where $\exp$ is related to $g$:
\be{401a}
c_\tau(t) = \exp_{\gamma_{x,\xi}(t)} \left( \tau v(t)  \right), \quad 
v(t) = \exp_{\gamma_{x,\xi}(t)}^{-1}\left(\hat \gamma_{x,\xi}(t)\right),
\ee
where $0\le t\le 1$, $0\le\tau\le1$. Then $c_0= \gamma_{x,\xi}$, $c_1= \hat\gamma_{x,\xi}$. Set $g^\tau = g + \tau(\hat g_1-g)$. Let 
\be{401b}
E(\tau) = \int_0^1 \langle \dot c_\tau(t), \dot c_\tau(t) \rangle_{g_\tau}\, \d t,
\ee
where, in local coordinates, $\langle \dot c_\tau(t), \dot c_\tau(t) \rangle_{g^\tau} = g^\tau_{ij} (c_\tau ) \dot c_\tau^i \dot c_\tau^j$. Apply Taylor's formula 
\[
E(1) = E(0) +E'(0) +\int_0^1 (1-\tau) E''(\tau)\,\d \tau
\]
to get
\be{402}
E'(0) = -\int_0^1 (1-\tau) E''(\tau)\, \d \tau
\ee
because $E(0)=E(1)= l_m^2(x,\xi)$. Write
\[
\psi(\tau,s) = \int_0^1 \langle \dot c_\tau(t), \dot c_\tau(t) \rangle_{g^s}\, \d t 
= \int_0^1  g^s_{ij}(c_\tau) \dot c_\tau^j \dot c_\tau^i   \,\d t,
\]
where the second integrand is written in local coordinates. 
Then $E(\tau) = \psi(\tau,\tau)$. 
For $E'$ we get
\be{402a}
E'(\tau) = \psi_\tau(\tau,\tau) + \psi_s(\tau,\tau) .
\ee
Since $c_0 = \gamma_{x,\xi}$ is a critical curve for the energy functional, we get $\psi_s(0,0)=0$, therefore,
\[
E'(0) =  
\int_0^1 \langle \tilde f, \dot\gamma_{x,\xi}^2\rangle\, \d t,
\]
recall \r{ft}. Together with \r{402} this yields
\be{403}
I\tilde f(\gamma_{x,\xi}) =  -\int_0^1 (1-\tau) E''(\tau)\, \d \tau.
\ee
To estimate the r.h.s.\ above, note that
\be{403aa}
E''(\tau) = \psi_{\tau\tau}(\tau,\tau) +2\psi_{\tau s}(\tau,\tau) 
\ee
because $\psi_{ss}=0$. 
Note that 
\be{403a}
\bigg| \frac{\partial  c_\tau(t)}{\partial\tau}\bigg| +
\bigg| \frac{\partial \dot c_\tau(t)}{\partial\tau}\bigg|  \le C\big(  |v(t)| +  |\dot v(t)|  \big),
\ee
and  
\be{404}
\bigg| \frac{\partial^2  c_\tau(t)}{\partial\tau^2}\bigg| +
\bigg| \frac{\partial^2 \dot c_\tau(t)}{\partial\tau^2}\bigg|  \le C\big(  |v(t)| +  |\dot v(t)|  \big)^2.
\ee
We have that $|\exp_x^{-1}y| \le C|x-y|$ for $|x-y|\ll1$, where the norm is in any fixed coordinate chart, and $C>0$ depends on an upper bound of $g$ in $C^k$, $k\gg2$. All constants below will have the same property. This and \r{401a} imply
\be{405}
|v(t)| \le C \left|\hat \gamma_{x,\xi}(t) - \gamma_{x,\xi}(t) \right|\le C' \|\tilde f\|_{C^1}.
\ee
Since in fixed coordinates, $D_x \exp_x^{-1}y + D_y \exp_x^{-1}y=0$ when $x=y$, we have 
\[
\big|D_x \exp_x^{-1}y + D_y \exp_x^{-1}y\big|\le C|x-y|.
\]
This allows us to estimate $\dot v(t)$, see \r{401a}, as follows
\be{406}
|\dot v(t)| \le C \left(  \big|\hat \gamma_{x,\xi}(t) - \gamma_{x,\xi}(t) \big|  +  \big|\dot{\hat \gamma}_{x,\xi}(t) - \dot\gamma_{x,\xi}(t) \big| \right)  \le  C' \|\tilde f\|_{C^1}.
\ee
By \r{403}, \r{403aa}, \r{403a}, \r{404}, \r{405}, \r{406},
\[
|I\tilde f(\gamma_{x,\xi}) | \le C' \|\tilde f\|_{C^1}^2.
\]
This is the same estimate that was used in the linearization argument in \cite{SU-Duke, SU-rig} and goes back to \cite{E}, but now proven in this more general situation. Therefore, 
\be{I5}
\|I_{\alpha_j} \tilde f\|_{L^\infty} \le C\|\tilde f\|_{C^1}^2.
\ee
That  implies the same for $\|N_{\alpha_j} \tilde f\|_{L^\infty}$, see e.g., \cite{SU-rig}, therefore, 
\be{I6}
\|N_{\alpha_j}\tilde  f\|_{L^\infty} \le C\|\tilde f\|_{C^1}^2, \quad\forall j.
\ee

\subsection{  End of the proof of  Theorem~\ref{thm_1}.   }
We will use interpolation to estimate $\|N_{\alpha_j} \tilde f\|_{\tilde H^2(\tilde{M})}$ through some power of $\|N_{\alpha_j}\tilde  f\|_{L^\infty}$.  Since $N_{\alpha_j}$ is a \PDO\ of order $-1$ and $\tilde f$, extended as $0$ outside $M$, is smooth enough if $k\gg2$ by \r{I3}, we get
\be{I7}
\|N_{\alpha_j} \tilde f\|_{\tilde H^2(\tilde{M})}^2  \le \|N_{\alpha_j} \tilde f\|_{H^2(\tilde{M})}^2 
\le C\|\tilde f\|_{H^3(M)} \|N_{\alpha_j} \tilde f\|_{L^2(\tilde{M})} .
\ee
Combine \eqref{I6} and \eqref{I7} to get
\be{I8}
\|N_\alpha \tilde f\|_{\tilde H^2(\tilde{M})}^2 \le C \|\tilde f\|_{C^{3}} \|\tilde f\|_{C^1}^2 \le C' \|f\|_{C^5}^3,
\ee
by \eqref{I2a2}.

Since $I_{g_0,\D'}$ is s-injective, so is $N_\alpha$, related to $g_0$, by the support properties of $\alpha$.  
Now, since $N_\alpha$ by \r{thm3} is s-injective for $g=g_0$, we get from  \cite[Theorem~2]{SU-AJM} that $N_\alpha$ (the one related to $g$) is s-injective as well provided that
\be{I8a}
\|g-g_0\|_{C^k(M)}\le \eps_0
\ee
with some $k\gg1$ and $\eps_0>0$. Moreover, \r{401_1} is true. Assume now that both \r{I8a} and \r{thm3} are satisfied. The geodesics issued from $\supp\alpha$ form a complete set by the second statement of Proposition~\ref{pr_1}, therefore,  by \r{I8} and \r{401_1}, 
\be{I9}
\|\tilde f^s\|_{L^2(M)} \le C\|N_\alpha \tilde f\|_{\tilde H^2} \le C'\|f\|_{C^5}^{3/2}.
\ee
By \eqref{I4}, $\tilde f^s = f+f_2^s$, therefore,
\[
\|\tilde f^s\|_{L^2(M)}\ge \|f\|_{L^2(M)}- \|f_2^s\|_{L^2(M)} \ge \|f\|_{L^2(M)}- C\|f\|^2_{C^2}.
\]
Together with \eqref{I9}, this yields
\[
\|f\|_{L^2(M)}\le C\left( \|f\|^2_{C^2} +\|f\|_{C^5}^{3/2}\right) \le C' \|f\|_{C^5}^{3/2}
\]
because the $C^5$ norm of $f$ is uniformly bounded when $\eps\le 1$. 
We can use now interpolation estimates in $C^k$, see \cite{Tr}, and Sobolev embedding estimates to get $\|f\|_{C^5}\le C\|f\|_{C^0}^\mu\le C'\|f\|_{H^{n/2+1}}^\mu$ with any $\mu\in(0,1)$ as long as $k=k(\mu)\gg1$ in \eqref{I2c}. Next, interpolation estimates in Sobolev spaces imply $\|f\|_{H^{n/2+1}} \le C \|f\|_{L^2}^\mu$, so in the end, we get
\[
\|f\|_{L^2(M)}  \le C \|f\|_{L^2(M)}^{3\mu^2/2}.
\]
Choose $2/3<\mu^2<1$ with a corresponding $k\gg2$, to get $\|f\|_{L^2(M)}\ge 1/C$ if $f\not =0$. 
Choose $\eps\ll1$ to get a contradiction with  \eqref{I2c}. This proves Theorem~\ref{thm_2} for $\eps$ replaced by $\min(\eps_0,\eps)$. 

Now, $f=0$ implies $h=g$, therefore, $g$ and $\hat g$ are isometric. 

This concludes the proof of Theorem~\ref{thm_1}. 

\section{Proof of Theorem~\ref{thm_2} and Corollary~\ref{cor_1}} \label{sec_rem1}

\begin{proof}[Proof of Theorem~\ref{thm_2}]
As we mentioned in the Introduction, Theorem~\ref{thm_2} is a reformulation of \cite[Theorem~3]{SU-AJM}, as we show below. 
Notice first, that the s-injectivity of $I_{g,\D''}$ for some $\D''\subset \D'$ implies s-injectivity of $I_{g,\D'}$. One can always assume that $M$ is equipped with a finite analytic atlas. 
Note that the assumption that $\D'$ is open implies that the corresponding set of geodesic is open in the sense of \cite{SU-AJM} by Lemma~\ref{lemma_open}. 
Let $\D''\Subset \D'$ be such that $\D''$ is still complete and open.  It can be constructed as in Section~\ref{sec_4.2}.  As in  Section~\ref{sec_4.2} again, with $\D''$ and $\D'$ playing the roles of $\D'$ and $\D$ there, respectively, we construct $N_\alpha$ such that all geodesics through $\supp\alpha$ cover $\kappa_-^{-1}(\D'')$ and are contained in the interior of $\kappa_-^{-1}(\D')$. Then $N_\alpha$ is s-injective for each analytic $g_0\in\mathcal{G}$ by \cite[Theorem~1]{SU-AJM}. It is still s-injective under a small enough $C^k$, $k\gg2$, perturbation of $g\in C^k(M)$ by \cite[Theorem~2]{SU-AJM}. Note that \cite[Theorem~2]{SU-AJM} requires that the perturbation must be considered in $C^k(\tilde M)$ but one can use extensions of $g$ near a fixed $g_0$ so that their norms in $C^k(\tilde M)$ are bounded by $\|g\|_{C^k(M)}$ with a fixed $C$. Using the fact that $\D''\Subset \D'$, Lemma~\ref{lemma_open}, and the support properties of $\alpha$, we deduce that $I_{g,\D'}$ is s-injective for $g$ close enough in $C^k(M)$ to a fixed analytic $g_0\in\mathcal{G}$. We can therefore build $\mathcal{G}_s$ as a small enough  neighborhood of the analytic metrics in $\mathcal{G}$, that form a  dense set. 

This completes the proof of Theorem~\ref{thm_2}.
\end{proof}

Corollary~\ref{cor_1} is an immediate consequence of Theorem~\ref{thm_1} and Theorem~\ref{thm_2}.

\begin{proof}[Proof of Remark~\ref{rem_1}]
Condition (T) is not needed for \cite[Theorem~2]{SU-AJM}, see also \r{401_1} but it is used in the proof of Theorems~2 (s-injectivity for analytic metrics) and Theorem~3 (generic s-injectivity) there. Assume now that (T) in Theorem~\ref{thm_2} of this paper is replaced by the assumption that  $(\tilde M, \tilde g)$ is simple as in the remark. In the proof above, in this situation, we need to show that $N_\alpha$ is s-injective for a dense set of metrics in $\mathcal{G}$. Fix $g_0\in \mathcal{G}$. It can be extended as $\tilde g_0$ to $\tilde M$ so that $(\tilde M, \tilde g)$ is simple. Given $\eps>0$ we can find a real analytic $\tilde g$ so that $\|\tilde g - \tilde g_0\|_{C^k(\tilde M)}\le \eps$ for $k$ fixed. Then the ray transform $I$ related to $\tilde g$ over all maximal geodesics in $\tilde M$ is s-injective, see \cite{SU-rig, SU-AJM}. Let now $N_{g,\alpha}f=0$ in $M$, where $g=\tilde g|_M$, and the subscript $g$ in $N_{g,\alpha}$ indicates that this is the normal operator related to $g$. Then we get after integration by parts that $I_{g}f(\gamma)=0$ for all maximal geodesics in $M$. Let $\tilde f$ be the extension of $f$ as zero outside $M$. Then $I_{\tilde g}\tilde f(\gamma)=0$ for all maximal geodesics in $\tilde M$. Therefore, $\tilde f = d\tilde v$ in $\tilde M$, $\tilde v\in H_0^1(M)$. Since $\tilde f=0$ in $\tilde M\setminus M$, on can see that the same holds for $\tilde v$ as well, see also \cite{Sh-sib}. Therefore, $f$ is potential, thus $I_{g}$ is s-injective, and so is $N_{g,\alpha}$, see \cite[Lemma~2]{SU-AJM}.
\end{proof}

\end{document}